\documentclass{amsart}

\usepackage{amsmath, amssymb,verbatim,appendix}
\usepackage{enumerate}
\usepackage{cite}
\usepackage{tikz-cd}

\newtheorem{theorem}{Theorem}[section]
\newtheorem{lemma}[theorem]{Lemma}
\newtheorem{corollary}[theorem]{Corollary}
\newtheorem{proposition}[theorem]{Proposition}

\newtheorem*{claim}{Claim}
\newtheorem{fact}[theorem]{Fact}

\newtheorem*{theorem*}{Theorem}

\theoremstyle{definition}

\newtheorem{definition}[theorem]{Definition}

\newcommand{\abar}{\overline{a}}
\newcommand{\bbar}{\overline{b}}
\newcommand{\cbar}{\overline{c}}

\newcommand{\calC}{\mathcal{C}}

\newcommand{\calQ}{\mathcal{Q}}

\newcommand{\calU}{\mathcal{U}}

\newcommand{\wt}[1]{\widetilde{#1}}

\newcommand{\id}{\operatorname{id}}

\newcommand{\Aut}{\mathrm{Aut}}
\def\acl{\operatorname{acl}}
\def\dcl{\operatorname{dcl}}
\def\tp{\operatorname{tp}}

\newcommand{\bg}[3][\calC]{\Aut_{#2}(#3/#1)}

\def\forkindep{\mathrel{\raise0.2ex\hbox{\ooalign{\hidewidth$\vert$\hidewidth\cr\raise-0.9ex\hbox{$\smile$}}}}}

\title{Note on a differential algebra bound}

\author{L\'eo Jimenez}
\address{L\'eo Jimenez\\
The Ohio State University\\
Department of Mathematics\\
Math Tower\\
Columbus, OH 43210-1174\\
United States}
\email{jimenez.301@osu.edu}

\keywords{geometric stability, differentially closed fields, differential Galois theory, internality to the constants, orthogonality}
\subjclass[2020]{03C45, 03C98, 12H05, 12L12}

\date{\today}

\begin{document}

\begin{abstract}
    In a recent article, Freitag, Moosa and the author showed that in differentially closed fields of characteristic zero, if two types are nonorthogonal, then their $n+3$ and $m+3$ Morley powers are not weakly orthogonal, where $n$ and $m$ are their respective Lascar ranks. 
    
    In this short note, we prove that the bound is tight: there are such types with weakly orthogonal $n+2$ and $m+2$ Morley powers. The types in question were constructed by Freitag and Moosa as examples of types with degree of nonminimality 2. As interesting as our result are our methods: we rely mostly on Galois theory and some descent argument for types, combined with the failure of the inverse Galois problem over constant parameters.
\end{abstract}

\maketitle

\section{Introduction}

In \cite{freitag2023differential}, Freitag, Moosa and the author proved the following theorem about differentially closed fields of characteristic zero:

\begin{theorem*}[Freitag-J.-Moosa]
    Let $p,q \in S(F)$ be two finite rank stationary types over some differential field $F$, and suppose that $p $ and $q$ are not orthogonal. Then their Morley powers $p^{(U(p)+3)} $ and $q^{(U(q)+3)}$ are not weakly orthogonal.
\end{theorem*}

This result can be interpreted, in a more intuitive manner, as follows: if there is a differential-algebraic relation between tuples $\abar \models p^{(k)}$ and $\bbar \models q^{(l)}$, for some $k$ and $l$, then there is a differential-algebraic relation between tuples of length $U(p)+3$ and $U(q)+3$.

The goal of this note is to prove that this bound is tight:

\begin{theorem*}[In $\mathrm{DCF}_0$]
    For any $n,m \in \mathbb{N}$, there is a differential field $F$ and types $p,q \in S(F)$ with $U(p) = n$ and $U(q) = m$, such that $p \not\perp q$ but $p^{(n+2)} \perp^w q^{(m+2)}$.
\end{theorem*}

The types used to obtain this bound are given by the expected construction: the projectivisations of generic linear differential equations, considered by Freitag and Moosa in \cite[Subsection 4.2]{freitag2023bounding}, and generalizing an example of Jin and Moosa \cite[Example 4.2]{jin2020internality}. They constructed these types to witness that in $\mathrm{DCF}_0$, the degree of nonminimality can be 2. This has since then been proved to be the largest possible value by Freitag, Jaoui and Moosa in \cite{freitag2023degree}.

The main feature of such a type $p$ is that it is internal to the field $\calC$ of constants, and that the action of its binding group $\bg{}{p}$ is definably isomorphic to the natural action of $\mathrm{PGL}_{n+1}(\calC)$ on $\mathbb{P}^n(\calC)$. Because this action is generically $n+2$ transitive, meaning that the diagonal action of $\mathrm{PGL}_{n+1}$ on $(\mathbb{P}^n)^{n+2}$ has an orbit of maximal rank, we obtain that $p^{(n+2)}$ is weakly $\calC$-orthogonal. Pick another such type $q$ with binding group $\mathrm{PGL}_m$, we similarly have that $q^{(m+2)}$ is weakly $\calC$-orthogonal. 

Once given these types, the proof relies exclusively on proving a preservation of binding group result. More precisely, we consider the binding group of $p$, not over the constants, but over the family of type-definable sets consisting of the constants and $q$, denoted $\bg[q, \calC]{}{p}$. We show that this binding group is still, in some cases, definably isomorphic to $\mathrm{PGL}_{n+1}$. From this it is not hard to prove, and this is Lemma \ref{lem: triv-quo-bound} below, that $p^{(n+2)}$ and $q^{(m+2)}$ are weakly orthogonal.

In the $n \neq m$ case, this preservation result is easily seen to follow from simplicity of $\mathrm{PGL}_n$ and $\mathrm{PGL}_m$ and using the Galois-theoretic fact that the quotients $\bg{}{p}/\bg[q,\calC]{F}{p}$ and $\bg{}{q}/\bg[p,\calC]{}{q}$ are definably isomorphic. 

The $n = m$ case is more subtle. The core of our argument is to show, using uniform internality techniques from \cite{jaoui2023relative}, the following descent result: if $p^{(n+2)}$ and $q^{(n+2)}$ are not weakly orthogonal, then their binding groups are definably isomorphic to the binding group of a type over constant parameters. But this group is also definably isogenous to $\mathrm{PGL}_{n+1}(\calC)$, which is well-known to not be the binding group of any type over constant parameters.

The organization of this article is as follows: Section 2 contains model-theoretic preliminaries, and Section 3 contains the proof of the Theorem. 

\section{Model-theoretic generalities}

In this section, we consider an $\omega$-stable theory $T$, eliminating imaginaries, and work in a large saturated model $\calU \models T$. We will denote $p(\calU)$ the set of realizations of some type $p$ in $\calU$.

We fix some parameter set $A$ and an $A$-definable subset $\calC$, we will consider internality with respect to $\calC$. 

We assume familiarity with geometric stability theory, and in particular binding groups. A good reference is \cite{pillay1996geometric}. We briefly recall some facts about internality and binding groups, more can be found in \cite[Chapter 7, Section 4]{pillay1996geometric}. 

Let $\calQ$ be a family of type definable sets, each defined over a subset of $A$.

\begin{definition}
    A type $p \in S(A)$ is said to be (resp. almost) $\calQ$-internal if it is stationary and there are some $B \supset A$, some $a \models p$, some $c_1, \cdots , c_n$ realizing some partial type in $\calQ$, such that:
    \begin{itemize}
        \item $a \forkindep_A B$,
        \item $a \in \dcl(c_1, \cdots , c_n, B)$ (resp. $\acl$).
    \end{itemize}
\end{definition}

This of course covers the case where $\calQ = \{ \calC \}$, and we recover the more common notion of internality to the definable set $\calC$. The binding group of $p$ over $\calQ$ is, abstractly, the group:
\[\bg[\calQ]{A}{p}  := \left\{ \sigma\vert_{p(\calU)} : \sigma \in \bg[\calQ]{A}{\calU} \right\}\]
\noindent i.e. restrictions to $p$ of automorphisms of $\calU$ fixing $A$ and $\calQ$ pointwise. The following classic theorem (see \cite[Theorem 7.4.8]{pillay1996geometric}) is key:

\begin{fact}
    If $p \in S(A)$ is $\calC$-internal, then the $\bg[\calQ]{A}{p}$, and its action of $p$, are isomorphic to an $A$-definable group action. Moreover, the $A$-definable group $\bg[\calQ]{A}{p}$ is definably isomorphic to a group definable in $\calQ^{\mathrm{eq}}$.
\end{fact}

We also call that definable group the binding group of $p$ over $\calQ$. Remark that in a general stable context, we only get type-definability of the binding group: definability is a consequence of $\omega$-stability.

We will consider two $\calC$-internal types $p,q \in S(A)$ and will attempt to understand forking between realizations of $p$ and $q$. Let us briefly recall the various orthogonality notions that we will use:

\begin{definition}
    Let $p,q \in S(A)$ be two types. They are:
    \begin{itemize}
        \item weakly orthogonal if for any $b \models p$ and $c \models q$, we have $b \forkindep_A c$. We write $p \perp^w q$.
        \item orthogonal if they are stationary, and for any $D \supset A$, the types $p\vert_D$ and $q \vert_D$ are weakly orthogonal. We write $p \perp q$.
    \end{itemize}

    The type $p$ is:
    \begin{itemize}
        \item  weakly $\calC$-orthogonal if for any tuple $\cbar$ of elements of $\calC$ and $b \models p$, we have $a \forkindep_A \cbar$.
        \item $\calC$-orthogonal if it is stationary and for all $D \supset A$, the non-forking extension $p\vert_D$ is weakly $\calC$-orthogonal.
    \end{itemize}

\end{definition}

We have the following well-known fact (see \cite[Chapter 1, Lemma 4.3.1]{pillay1996geometric}):

\begin{fact}
    Let $p,q \in S(A)$ be two stationary types. If $p \not\perp q$, then $p^{(n)} \not\perp^w q^{(m)}$ for some $n ,m \in \mathbb{N}$.
\end{fact}

The theorem of Freitag, Moosa and the author mentioned in the introduction is giving, in the case of $\mathrm{DCF}_0$, a precise upper bound on $n$ and $m$, and our goal here is to show that it is the best possible.

To do so, we will use the binding groups $\bg{A}{p}$ and $\bg{A}{q}$, but also the binding groups $\bg[q, \calC]{A}{p}$ and $\bg[p,\calC]{A}{q}$ (i.e. binding groups over the families $\{ p, \calC\}$ and $\{ q , \calC \}$). 

We first recall how weak $\calC$-orthogonality relates to the binding group action:

\begin{fact}\label{fact: bing-trans-weak-orth}
    Let $p \in S(A)$ be a $\calC$-internal type. Then $p$ is weakly $\calC$-orthogonal if and only if $\bg{A}{p}$ acts transitively on $p$.
\end{fact}

We will also need the following (see \cite[Section 2]{eagles2024splitting} for a proof):

\begin{lemma}\label{lem: bind over family normal}
    The group $\bg[q, \calC]{A}{p}$ is $A$-definably isomorphic to an $A$-definable normal subgroup of $\bg{A}{p}$.
\end{lemma}

\noindent as well as the following consequence of the definable Goursat's lemma (see \cite[Lemma 2.10]{eagles2024splitting}):

\begin{lemma}\label{lem: binding iso quotients}
    Consider some $\calC$-internal types $p,q \in S(A)$. Then the two quotients $\bg{A}{p}/\bg[q,\calC]{A}{p}$ and $\bg{A}{q}/\bg[p,\calC]{A}{q}$ are $A$-definably isomorphic.
\end{lemma}

Remark that the previous two results are, in \cite{eagles2024splitting}, given for types over an algebraically closed set of parameters. However, an examination of the proofs shows that stationarity of these types suffices.

In the case where these quotients are trivial, we obtain:

\begin{lemma}\label{lem: triv-quo-bound}
    Suppose that $p,q \in S(A)$ are $\calC$-internal and that $\bg{A}{p} = \bg[q,\calC]{A}{p}$ and $\bg{A}{q} = \bg[p,\calC]{A}{q}$. If $n,m$ are integers such that $p^{(n)}$ and $q^{(m)}$ are weakly $\calC$-orthogonal, then $p^{(n)} \perp^w q^{(m)}$.
\end{lemma}

\begin{proof}
    Let $\abar_1, \abar_2 \models p^{(n)}$ and $\bbar_1, \bbar_2 \models q^{(m)}$. Because $p^{(n)}$ and $q^{(m)}$ are weakly $\calC$-orthogonal, by Fact \ref{fact: bing-trans-weak-orth} there are $\sigma \in \bg{A}{p}$ and $\tau \in \bg{A}{q}$ such that $\sigma(\abar_1) = \abar_2$ and $\tau(\bbar_1) = \bbar_2$. 

    By our assumptions, these can be extended to $\wt{\sigma} \in \Aut_A(\calU)$ fixing $q(\calC)$ pointwise and $\wt{\tau} \in \Aut_A(\calU)$ fixing $p(\calU)$ pointwise. Then $\wt{\tau} \circ \wt{\sigma} (\abar_1, \bbar_1) = \wt{\tau} (\abar_2, \bbar_1) = (\abar_2, \bbar_2)$. 

    In other words, we have obtained that $\abar_1 \bbar_1 \equiv_A \abar_2 \bbar_2$ for any $\abar_1, \abar_2 \models p^{(n)}$ and $\bbar_1, \bbar_2 \models q^{(m)}$, which implies $p^{(n)} \perp^w q^{(m)}$, as we could have picked $\abar_1$ and $\bbar_1$ independent over $A$.
\end{proof}

One immediate way to use this is to work with two groups having no non-trivial isomorphic quotients. We will only consider a special case, but we record the general property we need, as we suspect it could have other uses.

\begin{definition}
    Two $A$-definable groups $G_1$ and $G_2$ are \emph{immiscible} if there does not exist $A$-definable subgroups $H_i \lneq G_i$ such that $G_1/H_1$ and $G_2/H_2$ are $A$-definably isomorphic. 
\end{definition}

We get:

\begin{corollary}\label{cor: incomp-bound}
    Suppose that $p,q \in S(A)$ are $\calC$-internal and $\bg{A}{p}$ and $\bg{A}{q}$ are immiscible. Let $n,m$ be integers such that $p^{(n)}$ and $q^{(m)}$ are weakly $\calC$-orthogonal. Then $p^{(n)} \perp^w q^{(m)}$.
\end{corollary}

\begin{proof}
    From the immiscibility assumption and Lemma \ref{lem: binding iso quotients}, we obtain $\bg{A}{p} = \bg[q,\calC]{A}{p}$ and $\bg{A}{q} = \bg[p,\calC]{q}{A}$. Lemma \ref{lem: triv-quo-bound} allows us to conclude.
\end{proof}

We will use that weak orthogonality to the constants is preserved under non-forking extension over parameters orthogonal to the constants. The following well-known fact can be obtained by a straightforward forking computation:

\begin{fact}\label{fact: nf-over-orth-w-orth}
    Let $p \in S(A)$ and let $d$ be some tuple such that $\tp(d/A)$ is orthogonal to $\calC$. If $p$ is weakly orthogonal (resp. orthogonal) to $\calC$, then $p\vert_{Ad}$ is weakly orthogonal (resp. orthogonal) to $\calC$.
\end{fact}

We will need to understand how binding groups change under non-forking extensions. The following lemma can be considered a finer version of Fact \ref{fact: nf-over-orth-w-orth}:

\begin{lemma}\label{lem: ortho-nf-ext}
    Let $p \in S(A)$ be a $\calC$-internal type, and let $d$ be some tuple. Then $p\vert_{Ad}$ is $\calC$-internal, and $\bg{Ad}{p\vert_{Ad}}$ injects definably into $\bg{A}{p}$. If $\tp(d/A)$ is stationary and orthogonal to $\calC$, the binding groups are definably isomorphic.
\end{lemma}

\begin{proof}
    The first part is proven in \cite[Fact 2.5]{eagles2024splitting}: pick a Morley sequence $\bbar = b_1 , \cdots , b_n$ in $p\vert_{Ad}$ that is long enough to be a fundamental system of both $p$ and $p\vert_{Ad}$. The map $\iota$ given by sending $\sigma \in \bg{A}{p\vert_{Ad}}$ to the unique $\iota(\sigma) $ such that $\iota(\sigma)(\bbar) = \sigma(\bbar)$ is a definable injection. 

    If $\tp(d/A)$ is stationary and orthogonal to $\calC$, we show that this map is surjective. So let $\tau \in \bg{A}{p}$, and pick any extension $\wt{\tau} \in \Aut_A(\calU)$. 

    \begin{claim}
        The tuple $d$ is independent of every tuple consisting of realizations of $p$ and elements of $\calC$.
    \end{claim}

    \begin{proof}[Proof of claim]
        Suppose that $\bbar$ is a a tuple of realizations of $p$ and elements of $\calC$, and that $d$ forks with $\bbar$ over $A$. By $\calC$-internality, there is some $B \supset A$ and a tuple $\cbar$ of constants such that $\bbar \forkindep_A B$ and $\bbar$ and $\cbar$ are interalgebraic over $B$, and we may also assume that $d \forkindep_A B$. We see that $d$ forks with $\bbar$ over $B$, and by interalgebraicity, that it forks with $\cbar$ over $B$, which contradicts orthogonality.
    \end{proof}
    
    In particular, as $\tp(d/A)$ is stationary, this means that any two realizations of $\tp(d/A)$ have the same type over $A \cup p(\calU) \cup \calC$. By stable embeddedness, this implies that for any two realisations, there is some automorphism in $\Aut_A(\calU/p , \calC)$ taking one to the other.

    Applying this to get some $\wt{\sigma} \in \Aut_A(\calU/p, \calC)$ such that $\wt{\sigma}(\wt{\tau}(d)) = d$, we get that $\wt{\sigma} \circ \wt{\tau}$ fixes $\calC$ pointwise, fixes $d$, and its restriction to $p(\calU)$ is $\tau$. Because it fixes $d$, it also restricts to $p\vert_{Ad}(\calU)$, and this restriction has image $\tau$ under the map $\iota$, proving surjectivity of $\iota$. 
     
\end{proof}

We will understand uniformly definable families of types as definable fibrations:

\begin{definition}
    If $p \in S(A)$, a \emph{definable fibration} $f : p \rightarrow f(p)$ is a partial $A$-definable map with domain containing $p(\calU)$ such that every fiber $\tp(a/f(a)A)$, for any $a \models p$, is stationary. The type $f(p)$ is the unique type such that $f(p)(\calU) = f(p(\calU))$.
\end{definition}
  
 We will need to apply to fibrations some of the uniform internality machinery from \cite{jaoui2023relative} (see also \cite[Lemma 5.7]{jimenez2019groupoids}):

\begin{definition}\label{def: unif-orth}
    Let $p \in S(A)$ be a stationary type and $\pi : p \rightarrow \pi(p)$ be a definable fibration. Then $\pi$ is \emph{uniformly $\calC$-internal} if there is $a \models p$ and some tuple $e$ such that:
    \begin{itemize}
        \item $a \forkindep_A e$
        \item there are $c_1, \cdots , c_n \in \calC$ such that $a \in \dcl(e,c_1, \cdots , c_n, \pi(a),A)$.
    \end{itemize}
    It is \emph{almost uniformly $\calC$-internal} if we replace $\dcl$ by $\acl$ in the second line.
\end{definition}

This implies, in particular, that the fibers of $\pi$ are $\calC$-internal, but is strictly stronger: if we only assume that the fibers are $\calC$-internal, then $e$ depends on $\pi(a)$, and we can only ask for $a \forkindep_{\pi(a)A} e$. 

A strong way uniform internality can happen is when the map $\pi$ is equivalent to the projection map of a product:

\begin{definition}
     Let $p \in S(A)$ be a stationary type and $\pi : p \rightarrow \pi(p)$ be a definable fibration. Then $\pi$ is \emph{split} if there is some $\calC$-internal type $r \in S(A)$, some $a \models p$ and some $b \models r$ such that $b \in \dcl(aA)$, the tuples $\pi(a)$ and $b$ are independent over $A$, and $a \in \dcl(\pi(a), b)$. 

     Replacing all $\dcl$ by $\acl$, we obtain the notion of \emph{almost splitting}.
\end{definition}

The only result from \cite{jaoui2023relative} that we will need is Proposition 3.14 (again, this was stated under the assumption that $A = \acl(A)$, but only stationarity is needed):

\begin{proposition}\label{prop: unif-int-iff-split}
    Let $p \in S(A)$ be stationary and $\pi : p \rightarrow \pi(p)$ be an $A$-definable fibration. Suppose that $\pi(p)$ is $\calC$-orthogonal. Then $\pi$ is uniformly almost $\calC$-internal if and only if it is almost split.
\end{proposition}

Note that the orthogonality assumption is required: see \cite[Theorem 4.4]{jin2020internality} and the preceding discussion for an example of a uniformly internal fibration that does not almost split (non-splitting follows from \cite[Theorem 4.5]{eagles2024splitting}).

\section{Maximizing weak orthogonality}

We now let $T = \mathrm{DCF}_0$, fix some large saturated $\calU \models \mathrm{DCF}_0$ and let and $\calC$ be the field of constants. This section will assume some familiarity with the model theory of differentially closed fields of characteristic zero, see \cite[Chapter 2]{marker2017model} for an introduction to this subject. Our examples will come from a construction of Freitag and Moosa in \cite{freitag2023bounding}, itself a generalisation of \cite[Example 4.2]{jin2020internality}. We give a quick presentation of their construction, and direct the reader to \cite[Subsection 4.2]{freitag2023bounding} for more details. 

Let $B = (b_{i,j})_{1 \leq i,j \leq n}$ be an $n \times n$ matrix of differentially transcendental, differentially independent elements of $\calU$. Write $\bbar = (b_{i,j})$ for the tuple formed by the $b_{i,j}$ and consider $F = \mathbb{Q}\langle \bbar \rangle$, the differential field the $b_{i,j}$ generate over $\mathbb{Q}$. Let $q \in S(F)$ be the generic type of $\delta (X) = BX$, it is stationary and $\calC$-internal.

Freitag and Moosa show that the binding group action of $\bg{F}{q}$ is definably isomorphic to the action of $\mathrm{GL}_n(\calC)$ on $\calC^n \setminus \{ 0 \}$. They then projectivise to obtain another stationary type $p \in S(F)$ such that the action of its binding group is definably isomorphic to the action of $\mathrm{PGL}_{n}(\calC)$ on $\mathbb{P}^{n-1}(\calC)$. In particular, because the degree of generic transitivity of this action is $n+2$ (i.e. the diagonal action of $\bg{F}{p}$ on $p^{(n+2)}$ is transitive, but it is not transitive on $p^{(n+3)}$), Fact \ref{fact: bing-trans-weak-orth} implies that $p^{(n+2)}$ is weakly orthogonal to the constants, but $p^{(n+3)}$ is not. 

In the proof we will constantly use, without mentioning it, the fact that the type of any tuple of differentially transcendental, independent elements is orthogonal to $\calC$. Indeed, any such type has $U$-rank $\omega^k$ for some $k \in \mathbb{N}$, and therefore must be orthogonal to any finite rank definable set by Lascar's inequalities. 

\begin{theorem}
    For any $n,m \in \mathbb{N}$, there are a differential field $L$ and types $p,q \in S(L)$ of $U$-ranks $n$ and $m$ such that $p \not\perp q$ and $p^{(n+2)} \perp^w q^{(m+2)}$.
\end{theorem}

\begin{proof}
    We have to deal with the case $n=m$ separately, so first assume that $n \neq m$.

    By the previous discussion, there are differential fields $F_n,F_m < \calU$ and some $\calC$-internal types $p_n,p_m \in S(F)$ such that the actions of $\bg{F}{p_i}$ on $p_i$ are definably isomorphic to the natural actions of $\mathrm{PGL}_{i+1}(\calC)$ on $\mathbb{P}^{i}(\calC)$ for $i = n ,m$ (in particular $U(p_i) = i$). 
    
    The fields $F_i$ are generated, as differential fields, by differentially independent differentially transcendental tuples $\bbar_n$ and $\bbar_m$. Moreover, we may assume that $\bbar_n \forkindep \bbar_m$, which implies that $\tp(\bbar_n / \bbar_m)$ is orthogonal to the constants (and symmetrically so is $\tp(\bbar_m/\bbar_n)$) by Fact \ref{fact: nf-over-orth-w-orth}. 
    
    Let $L$ be the differential field generated by $F_n$ and $F_m$. By the previous discussion and Fact \ref{fact: nf-over-orth-w-orth}, we have that $\left(p_i\vert_L\right)^{(i+2)}$ is weakly orthogonal to the constants. By Lemma \ref{lem: ortho-nf-ext}, the binding group of $p_i\vert_L$ is definably isomorphic to $\mathrm{PGL}_{i+1}(\calC)$. Since $n \neq m$, the binding groups of $p_n\vert_L$ and $p_m\vert_L$ are definably simple and non-isomorphic, thus immiscible. Therefore $\left( p_n\vert_L\right)^{(n+2)} \perp^w \left(p_m\vert_L \right)^{(m+2)}$. As both $p_n$
     and $p_m$ are $\calC$-internal, so are $p_n\vert_L$ and $p_m\vert_L$, and in particular $p_n\vert_L \not\perp p_m\vert_L$.

    \medskip

    Now assume that $n=m$. Of course, the previous method fails as the binding groups are now definably isomorphic. Here, we will take inspiration from the proof of \cite[Proposition 4.9]{jaoui2023relative}. 

    Again, we let $\bbar$ be the tuple of differentially independent differential transcendentals of the discussion preceding the theorem and $F$ be the differential field generated by $\bbar$. We also consider the type $p \in S(F)$ constructed by Freitag and Moosa, with binding group action definably isomorphic to $\mathrm{PGL}_{n+1}(\calC)$ acting on $\mathbb{P}^n(\calC)$.

    Fix some realization $a \models p$ and consider $s = \tp(a,\bbar)$, as well as the projection $\pi : s \rightarrow \pi(s)$ on the $\bbar$ coordinate (so $\pi(s)$ is the type of $n^2$ independent differential transcendentals). Note that $s$ is stationary as both $\tp(\bbar)$ and $\tp(a/\mathbb{Q}\langle \bbar \rangle)$ are. Pick two independent realizations $\bbar_1, \bbar_2$ of $\tp(\bbar) = \pi(s)$, and some $a_1, a_2$ with $a_i,\bbar_i \models s$. We let $p_i = \tp(a_i/ \mathbb{Q}\langle \bbar_i \rangle)$.

    We want to apply Lemma \ref{lem: binding iso quotients}, and thus need types over the same parameters. We consider the non-forking extensions $p_1\vert_{\bbar_2}$ and $p_2\vert_{\bbar_1}$, and let $L = \mathbb{Q}\langle \bbar_1 , \bbar_2 \rangle$. As previously, the types $\tp(\bbar_1/\bbar_2)$ and $\tp(\bbar_2/\bbar_1)$ are orthogonal to $\calC$, therefore $\bg{L}{p_1\vert_{L}}$ and $\bg{L}{p_2 \vert_{L}}$ are both isomorphic to $\mathrm{PGL}_{n+1}(\calC)$ by Lemma \ref{lem: ortho-nf-ext}, and $\left( p_i\vert_L \right)^{n+2}$ is weakly $\calC$-orthogonal for $i = 1,2$ by Fact \ref{fact: nf-over-orth-w-orth}.

    By Lemma \ref{lem: binding iso quotients} and simplicity of $\mathrm{PGL}_{n+1}$, we have either $\bg[p_2 \vert_L, \calC]{L}{p_1\vert_{L}} = \bg{L}{p_1 \vert_L}$ (and the symmetric equality for $p_2\vert_L$) or $\bg[p_2 \vert_L, \calC]{L}{p_1\vert_{L}} = \{ \id \}$. In the first case, we conclude by Lemma \ref{lem: triv-quo-bound} that $\left( p_1\vert_L \right)^{(n+2)} \perp^w \left( p_2\vert_L\right)^{(n+2)}$.

    So we just need to rule out the second case, which we assume for a contradiction. Note that this implies that $p_1\vert_L(\calU) \subset \dcl(p_2\vert_L(\calU), \calC, L)$. Fixing a fundamental system of solutions $\abar$ for $p_2\vert_L$ (i.e. a tuple of realizations of $p_2\vert_L$ such that $p_2\vert_L \subset \dcl(\abar, L ,\calC)$, which exists by $\calC$-internality of $p_2$), we obtain $p_1\vert_L(\calU) \subset \dcl(\abar, L, \calC)$. Thus the map $\pi : s \rightarrow \pi(s)$ is uniformly $\calC$-internal. As $\pi(s)$ is orthogonal to $\calC$, we have that $\pi$ almost splits by Proposition \ref{prop: unif-int-iff-split}. Let $r \in S(\mathbb{Q})$ be the $\calC$-internal type witnessing it. It is weakly $\calC$-orthogonal because $s$ is.

    Consider again the original $(a,\bbar) \models s$, by splitting there is $d \models r$ such that $(a,\bbar)$ and $(d,\bbar)$ are interalgebraic over $\mathbb{Q}$ and $d \forkindep \bbar$. By Lemma \ref{lem: ortho-nf-ext}, the binding groups $\bg{\mathbb{Q}}{r}$ and $\bg{F}{\tp(d/F)}$ are definably isomorphic. By Lemma \ref{lem: binding iso quotients} (see also \cite[Corollary 2.13]{eagles2024splitting}), the binding groups $\bg{F}{\tp(a/F)}$ and $\bg{F}{\tp(d/F )}$ are isogenous, the former being definably isomorphic to $\mathrm{PGL}_{n+1}(\calC)$. Combining these facts, we get that $\bg{\mathbb{Q}}{r}$ is isogenous to $\mathrm{PGL}_{n+1}(\calC)$, and in particular is linear. But any $\calC$-internal, weakly $\calC$-orthogonal type over some constant parameters, if it has a linear binding group, actually has an abelian binding group (see for example the proof of Theorem 3.9 in \cite{freitag2022any}), a contradiction

\end{proof}

\bibliography{biblio}
\bibliographystyle{plain}

\end{document}